\documentclass[12pt]{article}
\usepackage{amssymb, amsmath}
\begin{document}
\renewcommand{\theequation}{\arabic{section}.\arabic{equation}}
\newtheorem{thm}{Theorem}[section]
\newtheorem{lem}{Lemma}[section]
\newtheorem{deff}{Definition}[section]
\newtheorem{pro}{Proposition}[section]
\newtheorem{cor}{Corollary}[section]
\newcommand{\n}{\nonumber}
\newcommand{\tv}{\tilde{v}}
\newcommand{\tw}{\tilde{\omega}}
\renewcommand{\t}{\theta}
\newcommand{\w}{\omega}
\newcommand{\g}{\gamma}
\newcommand{\e}{\varepsilon}
\renewcommand{\a}{\alpha}
\renewcommand{\l}{\lambda}
\newcommand{\vare}{\varepsilon}
\newcommand{\s}{\sigma}
\renewcommand{\o}{\omega}
\renewcommand{\O}{\Omega}
\newcommand{\bb}{\begin{equation}}
\newcommand{\ee}{\end{equation}}
\newcommand{\bq}{\begin{eqnarray}}
\newcommand{\eq}{\end{eqnarray}}
\newcommand{\bqn}{\begin{eqnarray*}}
\newcommand{\eqn}{\end{eqnarray*}}
\title{On the blow-up problem and new \textit{a priori} estimates for the 3D Euler  and the Navier-Stokes equations}
\author{Dongho Chae\thanks{This research was supported partially
by  KRF Grant(MOEHRD, Basic Research Promotion Fund); $(\dagger)$ permanent address.\newline
}\\
Department of Mathematics\\
University of Chicago\\
Chicago, IL 60637, USA\\
  e-mail: {\it chae@math.uchicago.edu}\\
and \\
 Sungkyunkwan University$^{(\dagger)}$\\
Suwon 440-746, Korea}
 \date{}
\maketitle
\begin{abstract}
We study  blow-up rates and the blow-up profiles of possible asymptotically self-similar singularities of the 3D Euler equations, where the  sense of convergence and self-similarity
are considered in various sense. We extend much further, in particular, the previous nonexistence results of  self-similar/asymptotically self-similar singularities obtained in \cite{cha1,cha2}. Some implications the notions for the 3D Navier-Stokes equations are also deduced.
Generalization of the self-similar transforms is also considered, and
by appropriate choice of the transform  we obtain new \textit{a priori} estimates
for the 3D Euler and the Navier-Stokes equations.
\end{abstract}
\section{Asymptotically self-similar singularities}
 \setcounter{equation}{0}

We are  concerned  on the following Euler equations for the
homogeneous incompressible fluid flows in $\Bbb R^3$.
 \[
\mathrm{ (E)}
 \left\{ \aligned
 &\frac{\partial v}{\partial t} +(v\cdot \nabla )v =-\nabla p ,
 \quad (x,t)\in {\Bbb R^3}\times (0, \infty) \\
 &\quad \textrm{div }\, v =0 , \quad (x,t)\in {\Bbb R^3}\times (0,
 \infty)\\
  &v(x,0)=v_0 (x), \quad x\in \Bbb R^3
  \endaligned
  \right.
  \]
where $v=(v_1, v_2, v_3 )$, $v_j =v_j (x, t)$, $j=1,2,3$, is the
velocity of the flow, $p=p(x,t)$ is the scalar pressure, and $v_0
$ is the given initial velocity, satisfying div $v_0 =0$.
The system (E) is first modeled by  Euler in \cite{eul}.
The local well-posedness of the Euler equations in $H^m (\Bbb R^3)$, $m>5/2$, is established by Kato in \cite{kat}, which says that given $v_0 \in H^m (\Bbb R^3)$, there exists $T \in (0, \infty]$ such that
 there exists unique solution to (E), $v \in C([0, T);H^m (\Bbb R^3))$.
 The finite time blow-up problem of the local classical solution is known as one of the most
 important and difficult problems in partial differential equations(see e.g. \cite{maj,che,con1,con2,cha1} for graduate level texts and survey articles on the current status of the problem).
We say a local in time classical solution  $v\in C([0, T);H^m (\Bbb R^3))$ blows up at $T$ if
$ \lim\sup_{t\to T}\|v(t)\|_{H^m} =\infty
$
for all $m>5/2$.
The celebrated Beale-Kato-Majda criterion(\cite{bea}) states that the blow-up happens at $T$ if and only if
$$ \int_0 ^T \|\o (t )\|_{L^\infty} dt=\infty. $$
There are  studies of geometric  nature for the blow-up criterion(\cite{con3, con2, den}).  As another direction of studies of the blow-up problem mathematicians also consider various scenarios of singularities and study carefully their possibility of realization(see e.g. \cite{cor1,cor2, cha2,cha3} for some of those studies). One of the purposes in this paper, especially in this section, is to study more deeply the notions related to the scenarios of the self-similar singularities in the Euler equations, the preliminary studies of which are done in \cite{cha2, cha3}.
We recall that system (E) has scaling property that
  if $(v, p)$ is a
solution of the system (E), then for any $\lambda >0$ and $\alpha
\in \Bbb R $ the functions
 \bb
 \label{1.1}
  v^{\lambda, \alpha}(x,t)=\lambda ^\alpha v (\lambda x, \l^{\a +1}
  t),\quad p^{\l, \a}(x,t)=\l^{2\a}p(\l x, \l^{\a+1} t )
  \ee
  are also solutions of (E) with the initial data
  $ v^{\lambda, \alpha}_0(x)=\lambda ^\alpha v_0
   (\lambda x)$.
 In view of the scaling
  properties in (\ref{1.1}), a natural self-similar blowing up
  solution $v(x,t)$ of (E) should be of the form,
  \bq
  \label{1.2}
 v(x, t)&=&\frac{1}{(T-t)^{\frac{\a}{\a+1}}}
\bar{V}\left(\frac{x}{(T-t)^{\frac{1}{\a+1}}}\right)\\
\label{1.3}
p(x,t)&=&\frac{\a +1}{(T-t)^{\frac{2\a}{\a+1}}}
\bar{P}\left(\frac{x}{(T-t)^{\frac{1}{\a+1}}}\right)
 \eq
 for  $\a \neq -1$ and $t$  sufficiently
 close to $T$.
 Substituting (\ref{1.2})-(\ref{1.3}) into (E), we obtain the
following stationary system.
\bb\label{1.4}
\left\{ \aligned
  & \a \bar{V} +(y \cdot \nabla)\bar{V} + (\a+1)(\bar{V}\cdot \nabla )\bar{V} =-\nabla
  \bar{P},\\
 & \mathrm{div}\, \bar{V}=0,
  \endaligned \right.
\ee
the Navier-Stokes equations version of which has been studied extensively after Leray's pioneering paper(\cite{ler, nec, tsa, mil, cha3, hou}).
Existence of solution of the system (\ref{1.4}) is equivalent to the existence of solutions to the Euler equations of the form (\ref{1.2})-(\ref{1.3}), which blows up in a self-similar fashion.
Given $(\a,  p)\in (-1, \infty)\times  (0, \infty]$,  we say the blow-up is \emph{$\a-$asymptotically self-similar in the sense of $L^p$ }if there exists
$\bar{V}=\bar{V}_\a\in \dot{W}^{1,p}(\Bbb R^3 )$ such that the following convergence holds true.
$$
\lim_{t\to T} (T-t) \left\|\nabla v(\cdot, t)-\frac{1}{T-t} \nabla \bar{V}\left( \frac{\cdot}{(T-t)^{\frac{1}{\a+1}}}\right)\right\|_{L^\infty}=0
$$
if $p=\infty$, while
$$
\lim_{t\to T}(T-t)^{1-\frac{3}{(\a+1)p}}\left\|\o(\cdot ,t)-\frac{1}{(T-t)^{1-\frac{3}{(\a+1)p}}}\bar{ \O }\left(\frac{\cdot }{(T-t)^{\frac{1}{\a+1}}}\right)\right\|_{L^p}=0
$$
if $0<p<\infty$, where and hereafter we denote
$$\O=\mathrm{curl} \,V \quad \mbox{and}\quad \bar{ \O}=\mathrm{curl}\, \bar{V}.
$$
The above limit function $\bar{V}\in L^p(\Bbb R^3)$ with $ \bar{\O }\neq 0$ is called the \emph{blow-up profile}.
We observe that the self-similar blow-up given by  (\ref{1.2})-(\ref{1.3}) is trivial case of $\a-$asymptotic
self-similar blow-up with the blow-up profile given by the representing function $\bar{V}$.
We say a blow-up at $T$ is of \emph{type I}, if
 $$
 \lim\sup_{t\to T} (T-t)\|\nabla v (t)\|_{L^\infty}<\infty.
 $$
 If the blow-up is not of type I, we say it is of \emph{type II}.
For the use of terminology, type I and type II blow-ups, we followed the literatures on the studies of the blow-up problem in the semilinear heat equations(see e.g. \cite{mat,gig1,gig2}, and references therein).
The use of $\|\nabla v(t)\|_{L^\infty}$ rather than $\|v(t)\|_{L^\infty}$ in our definition of type I and II is motivated by  Beale-Kato-Majda's blow-up criterion.
\begin{thm}
Let $ m >5/2$, and $v\in C([0, T); H^m (\Bbb R^3 ))$ be a solution to  (E) with $v_0 \in H^m (\Bbb R^3)$, div $v_0=0$.
We set
\bb\label{1.5}
\lim\sup_{t\to T} (T-t)\|\nabla v (t)\|_{L^\infty}:=M(T).
\ee
Then, either $M(T)=0$ or $M(T)\geq 1$. The former case corresponds to non blow-up, and the latter case corresponds to the blow-up at $T$. Hence, the blow-up at $T$ is of type I if and only if $M(T) \geq 1$.
\end{thm}
{\bf Proof }
It suffices to show that $M(T)<1$ implies non blow-up at $T$, which, in turn, leads to $M(T)=0$, since $\|\nabla v(t)\|_{L^\infty}\in C([0, T])$ in this case.
 We suppose $M(T)<1$. Then, there exists
$t_0 \in (0, T)$ such that
$$\sup_{t_0 <t<T} (T-t)\|\nabla v (t)\|_{L^\infty}:=M_0 <1.$$
Taking curl of the evolution part of (E), we have
 the vorticity equation,
 $$
 \frac{\partial \o }{\partial t} +(v\cdot \nabla )\o =(\o
\cdot
 \nabla )v .
 $$
This, taking dot product with $\xi=\o/|\o|$, leads to
$$
 \frac{\partial |\o |}{\partial t} +(v\cdot \nabla )|\o |= (\xi\cdot \nabla)v \cdot \xi |\o|.
$$
Integrating this over $[t_0, t]$ along the particle trajectories $\{
X(a,t)\}$ defined by $v(x,t)$, we have
 \bb\label{1.6}
 |\o (X(a,t),t)|=|\o (X(a,t_0), t_0)|\exp \left[ \int_{t_0} ^t (\xi\cdot \nabla)v \cdot \xi
(X(a,s),s) ds \right],
  \ee
 from which we estimate
\bq\label{1.7}
\|\o (t)\|_{L^\infty}&\leq &\|\o(t_0) \|_{L^\infty} \exp\left[\int_{t_0}^t \|\nabla v (\tau)\|_{L^\infty} d\tau\right]\n\\
&<&\|\o (t_0) \|_{L^\infty}\exp\left[M_0 \int_{t_0} ^t
(T-\tau )^{-1} d\tau \right]\n\\
&=&\|\o (t_0)  \|_{L^\infty}\left(\frac{T-t_0}{T-t}\right)^{M_0}.
\eq
Since $M_0<1$, we have
$
\int_{t_0}^T\|\o (t)\|_{L^\infty} dt <\infty,
$
 and thanks to the Beale-Kato-Majda criterion there exists no blow-up at $T$, and  we can continue our classical solution beyond $T$.
 $\square$\\
\ \\
The following is our main theorem in this section.
 \begin{thm} Let a classical solution $v\in C([0, T); H^m (\Bbb R^3))$ with initial data $v_0 \in H^m (\Bbb R^3)\cap  \dot{W}^{1, p} (\Bbb R^3)$, div $v_0=0$, $\o_0 \neq 0$  blows up with type I. Let $M=M(T)$ be as in Theorem 1.1.
Suppose  $(\a,p)\in (-1, \infty)\times (0, \infty]$ satisfies
  \bb\label{1.8}
M <
\left| 1-\frac{3}{(\a+1)p}\right|.
\ee
 Then, there exists no $\a-$asymptotically self-similar blow-up at $t=T$ in the sense of $L^p$ if $\o_0 \in L^p (\Bbb R^3)$. Hence, for any type I blow-up and for any $\a \in (-1, \infty)$ there exists $ p_1\in  (0,\infty]$ such that  it is not $\a-$asymptotically self-similar  in the sense of $L^{p_1}$.
\end{thm}
{\bf Remark 1.1} We note that the case $p=\infty$ of the above theorem follows from Theorem 1.1, which states that there is no singularity at all at $t=T$ in this case.
The above theorem can be regarded an improvement of the main theorem  in \cite{cha3}, in the sense that we can consider the $L^p$ convergence only to exclude nontrivial blow-up profile $\bar{V}$, where $p$ depends on $M$. Moreover, we do not need to use the Besov space $\dot{B}^0_{\infty, 1}$ in the statement of the theorem, and the continuation principle of local solution in the Besov space in the proof.\\
\ \\
 \noindent{\bf Proof of Theorem 1.2}  We assume asymptotically self-similar blow-up happens at $T$.
Let us introduce  similarity variables defined by
$$
y=\frac{x}{(T-t)^{\frac{1}{\a+1}}}, \quad s=\frac{1}{\a+1} \log
\left( \frac{T}{T-t}\right),
$$
and transformation of the unknowns
$(v,p) \to (V, P)$ according to
 \bb\label{1.9}
 v(x,t)=\frac{1}{(T-t)^{\frac{\a}{\a+1}}} V(y,s ), \quad
 p(x,t)=\frac{1}{(T-t)^{\frac{2\a}{\a+1}}} P(y,s ).
 \ee
 Substituting $(v,p)$ into the $(E)$ we obtain the equivalent evolution equation for
 $(V,P)$,
 $$
 (E_1) \left\{ \aligned
  & V_s +\a V +(y \cdot \nabla)V + (\a+1)(V\cdot \nabla )V =-\nabla
  P,\\
 & \mathrm{div}\, V=0,\\
 & V(y,0)=V_0 (y)=T^{\frac{\a}{\a+1}} v_0 (T^{\frac{1}{\a}}y).
  \endaligned \right.
 $$
Then the assumption of asymptotically self-similar singularity at $T$ implies that there exists
 $\bar{V}=\bar{V}_\a\in \dot{W}^{1,p} (\Bbb R^3)$ such that
\bb\label{1.10}
 \lim_{s\to \infty}\|\O (\cdot, s)-\bar{\O}\|_{L^p}=0.
 \ee
 Now the  hypothesis (\ref{1.8}) implies that there exists  $t_0 \in (0, T)$ such that
 \bb\label{1.10a}
 \sup_{t_0 <t<T} (T-t)\|\nabla v (t)\|_{L^\infty} :=M_0 <
\left| 1-\frac{3}{(\a+1)p}\right|.
\ee
Taking $L^p(\Bbb R^3)$ norm of (\ref{1.6}), taking into account the following simple estimates,
$$
 -\|\nabla v(\cdot ,t)\|_{L^\infty}\leq (\xi\cdot \nabla)v \cdot \xi
(x,t) \leq \|\nabla v (\cdot ,t)\|_{L^\infty} \quad \forall (x,t)\in
\Bbb R^3\times [t_0 , T),
$$
we obtain, for all $p\in (0, \infty]$,
\bq\label{1.11}
 \lefteqn{\|\o (t_0 ) \|_{L^p} \exp \left[- \int_{t_0} ^t \|\nabla v
      (\cdot,s)\|_{L^\infty} ds \right]\leq \|\o (t)\|_{L^p}}\hspace{1.in}\n \\
      && \qquad \leq \|\o_0 \|_{L^p}\exp \left[ \int_{t_0} ^t \|\nabla v
      (\cdot,s)\|_{L^\infty} ds \right],
 \eq
where we use the fact that $a\mapsto X(a,t)$ is a volume preserving map.
From the fact
 $$
 \int_{t_0} ^t \|\nabla v
      (\cdot,s)\|_{L^\infty} ds\leq M_0 \int_{t_0}^t (T-\tau)^{-1}d\tau=-M_0 \log \left( \frac{T-t}{T-t_0}\right),
 $$
and
$$
\frac{\|\o (t)\|_{L^p}}{\|\o (t_0)\|_{L^p}}
= \left(\frac{T-t}{T-t_0}\right)^{ \frac{3}{(\a+1)p} -1} \frac{\|\O(s)\|_{L^p}}{\|\O (s_0 )\|_{L^p}},
$$
where we set
$$s_0 =\frac{1}{\a+1} \log \left(\frac{T}{T-t_0}\right),
$$
 we find that (\ref{1.1}) leads us to
 \bb\label{1.12}
 \left(\frac{T-t}{T-t_0}\right)^{M_0+1- \frac{3}{(\a+1)p}}
 \leq \frac{\|\O(s)\|_{L^p}}{\|\O (s_0)\|_{L^p}}
 \leq  \left(\frac{T-t}{T-t_0}\right)^{-M_0+1- \frac{3}{(\a+1)p}}
\ee
for all $p\in (0, \infty]$.
Passing $t\to T$, which is equivalent to $s\to \infty$ in (\ref{1.12}), we have
from (\ref{1.10})
\bb\label{1.12a}
\lim_{s\to \infty} \frac{\|\O(s)\|_{L^p}}{\|\O (s_0)\|_{L^p}}
=\frac{\|\bar{\O}\|_{L^p}}{\|\O (s_0 ) \|_{L^p}}\in (0, \infty).
\ee
By (\ref{1.10a}) $M_0+1- \frac{3}{(\a+1)p}<0$ or $-M_0+1- \frac{3}{(\a+1)p} >0$.
In the former case we have
\bb\label{1.13}
\lim_{t \to T} \left(\frac{T-t}{T-t_0}\right)^{M_0+1- \frac{3}{(\a+1)p}}=\infty,
\ee
while, in the latter case
\bb\label{1.14}
\lim_{t \to T} \left(\frac{T-t}{T-t_0}\right)^{-M_0 +1- \frac{3}{(\a+1)p}} =0.
\ee
Both of (\ref{1.13}) and (\ref{1.14}) contradicts with (\ref{1.12a}).
If the blow-up is of type I, and $M(T)<\infty$, then one can always choose
 $p_1\in (0, p_0)$ so small that (\ref{1.8}) is valid for $p=p_1$.  With such $p_1$ it is not  $\a-$asymptotically self-similar in $L^{p_1}$.
 $\square$\\
\ \\
 For the self-similar blowing-up solution of the form (\ref{1.2})-(\ref{1.3}) we observe that in order to
 be consistent with the energy conservation, $\|v(t)\|_{L^2}=\|v_0\|_{L^2}$ for all $t\in [0, T)$, we need to fix $\a=3/2$.  Since the self-similar blowing up solution corresponds to  a trivial convergence of the asymptotically self-similar blow-up,
the following is immediate from Theorem 1.2.
  \begin{cor}
Given $p\in (0, \infty]$, there exists no self-similar blow-up with the blow-up profile $V$ satisfying $\O \in L^p (\Bbb R^3)$ if
  \bb\label{1.15}
\|\nabla V\|_{L^\infty} <
\left| 1-\frac{6}{5p}\right|.
\ee
   \end{cor}
   {\bf Remark 1.2 } The above corollary implies that we can exclude self-similar singularity of the Euler equations only under the assumption of $\O \in L^p (\Bbb R^3)$ if $p$ satisfies the condition (\ref{1.15}). \\
   \ \\
   The following is, in turn, immediate from the above corollary, which is nothing but Theorem 1.1 in \cite{cha2}.
  \begin{cor}
  There exists no self-similar blow-up with the blow-up profile $V$ satisfying $\O \in L^p (\Bbb R^3)$ for all $p\in (0, p_0 )$ for some $p_0 >0$.
   \end{cor}
   The following theorem is concerned on the possibility of type II asymptotically self-similar
   singularity of the Euler equations, for which  the blow-up rate near the possible  blow-up time $T$ is
 \bb\label{type2}\|\nabla v(t)\|_{L^\infty}\sim \frac{1}{(T-t)^\g }, \qquad \g >1.
 \ee
\begin{thm}
Let $v\in C([0, T); H^m (\Bbb R^3))$, $m>5/2$, be local classical solution of the Euler equations. Suppose  there exists $\g >1$ and $R_1 >0$ such that the following convergence holds true.
\bb\label{1.16}
\lim_{t\to T}(T-t)^{(\a-\frac32)\frac{\g}{\a+1}} \left\|v(\cdot, t)-\frac{1}{(T-t)^{(\a-\frac32)\frac{\g}{\a+1}}}\bar{V}\left(\frac{\cdot}{(T-t)^{\frac{\g}{\a+1}}} \right) \right\|_{L^2 (B_{R_1})}=0,
\ee
where $B_{R_1}=\{ x\in \Bbb R^3 \, | \, |x|<R_1\}$. Then, the blow-up profile $\bar{V}\in L^2_{loc} (\Bbb R^3 )$ is a weak solution of the following stationary  Euler equations,
\bb
(\bar{V}\cdot \nabla)\bar{V}=-\nabla \bar{P}, \qquad \mathrm{div} \,\bar{V}=0.
\ee
\end{thm}
\noindent{\bf Proof }
We introduce a self-similar transform defined by
\bb\label{1.17}
 v(x,t)=\frac{1}{(T-t)^{\frac{\a \g}{\a+1}}}   V\left(y, s \right),\quad
 p(x,t)=\frac{1}{(T-t)^{\frac{2\a \g}{\a+1}}} P\left(y, s\right)
 \ee
 with
\bb\label{1.18} y=\frac{1}{(T-t)^{\frac{\g}{\a+1}}} x, \quad
 s=\frac{1}{(\g-1 )T^{\g-1}} \left[ \frac{T^{\g-1}}{(T-t)^{\g -1}} -1\right].
\ee
Substituting $(v,p)$ in (\ref{1.17})-(\ref{1.18}) into the $(E)$, we have
\bb\label{1.19}
  (E_2) \left\{ \aligned
  & -\frac{\g}{s(\g-1) +T^{1-\g}}\left[ \frac{\a}{\a+1} V +\frac{1}{\a+1}(y \cdot \nabla)V \right]=V_s+(V\cdot \nabla )V +\nabla P,\\
 & \qquad \mathrm{div}\, V=0,\\
 & V(y,0)=V_0(y)=v_0 (y).
  \endaligned \right.
 \ee
 The hypothesis (\ref{1.16}) is written as
\bb\label{1.20}
\lim_{s\to \infty}\|V(\cdot, s)-\bar{V}(\cdot )\|_{L^2 (B_{R(s)})}=0,\quad R(s)=\left[ (\g-1 )s +\frac{1}{T^{\g-1}} \right]^{\frac{\g}{(\a +1)(\g-1)}},
\ee
which implies that
 \bb\label{1.21}
 \lim_{s\to \infty} \|V(\cdot , s)-\bar{V}\|_{L^2 (B_R)}=0, \qquad \forall R>0,
 \ee
 where $V(y,s) $ is defined by (\ref{1.17}).
Similarly to \cite{hou, cha3}, we consider the scalar test function $\xi
\in C^1_0 (0,1)$ with $\int_0 ^1\xi (s)ds\neq 0$,  and the vector
test function $\phi =(\phi_1 , \phi_2, \phi_3 )\in C_0^1 (\Bbb
 R^3)$ with div $\phi=0$.

 We multiply the first equation of $(E_2)$,
 in the dot product, by $\xi
 (s-n)\phi (y)$, and integrate it over $\Bbb R^3\times [n, n+1]$,
 and then we  integrate by parts to obtain
 \bqn
 &&+\frac{\a}{\a+1}\int_0 ^{1}\int_{\Bbb R^3}g(s+n)\xi (s)  V(y, s+n)\cdot\phi(y) dyds \\ &&-\frac{1}{\a +1}\int_0 ^{1}\int_{\Bbb R^3}g(s+n)\xi (s)V(y, s+n)\cdot (y \cdot \nabla)\phi (y)dyds\\
 &&\qquad=\int_0^{1}\int_{\Bbb R^3} \xi _s(s) \phi(y)\cdot V(y,s+n)
 dyds\\
 &&+\int_0 ^{1}\int_{\Bbb R^3}\xi (s)\left[V(y,s+n)\cdot (V(y,s+n)\cdot \nabla )\phi (y)\right]  dyds=0,
 \eqn
 where we set
 $$g(s)=\frac{\g}{s(\g-1) +T^{1-\g}}.$$
  Passing to the limit $n\to \infty$ in this equation, using the facts
   $\int_0 ^1\xi _s(s)ds=0$, $\int_0 ^1\xi
  (s)ds\neq 0$, $V(\cdot, s+n)\to \bar{V}$ in $L^2_{\mathrm{loc}} (\Bbb R^3)$, and finally
   $g(s+n)\to 0$,
  we  find that $\bar{V}\in L^2_{\mathrm{loc}} (\Bbb R^3)$ satisfies
$$
\int_{\Bbb R^3}
\bar{V}\cdot (\bar{V}\cdot \nabla
 )\phi (y) dy=0
 $$
 for all vector test function $\phi \in C_0^1 (\Bbb
 R^3)$ with div $\phi=0$.
 On the other hand, we can pass $s\to \infty$ directly in the weak formulation of the
 second equation of $(E_2)$ to have
 $$
 \int_{\Bbb R^3} \bar{V}\cdot \nabla \psi (y)dy=0
 $$
for all scalar test function $\psi \in C^1_0 (\Bbb R^3)$. $\square$\\

\section{Generalized similarity transforms and new a priori estimates}
 \setcounter{equation}{0}

  Let us consider a classical solution to (E) $v\in C([0, T);H^m (\Bbb R^3))$, $m>5/2$, where
  we assume $T\in (0, \infty]$ is the maximal time of existence of the classical solution.
  Let $p(x,t)$ be the associated pressure.
  Let $\mu(\cdot)\in C^1([0, T))$ be a scalar function such that  $\mu(t)>0$ for all $t\in [0, T)$ and $\int_0 ^T \mu (t) dt=\infty$.
  We transform from
$(v,p)$ to $(V, P)$ according to the formula,
 \bq\label{2.1}
 v(x,t)&=&\mu (t)^{\frac{\a}{\a+1}} V\left(\mu(t)^{\frac{1}{\a+1}}x, \int_0 ^t \mu (\sigma )d\sigma \right),\\
 \label{2.2}
 p(x,t)&=&\mu (t)^{\frac{2\a}{\a+1}} P\left(\mu(t)^{\frac{1}{\a+1}}x, \int_0 ^t \mu (\sigma )d\sigma \right),
 \eq
 where $\a\in (-1, \infty)$ as previously.
This means that the  space-time variables are transformed from $(x,t)
\in \Bbb R^3 \times [0,T)$ into $(y,s)\in \Bbb R^3\times [0,
\infty)$ as follows:
\bb\label{2.3}
y=\mu(t)^{\frac{1}{\a+1}}x, \quad s=\int_0 ^t \mu (\sigma )d\sigma.
\ee
Substituting (\ref{2.1})-(\ref{2.3}) into the Euler equations, we obtain the equivalent equations satisfied by $(V,P)$
  $$
 (E_*) \left\{ \aligned
  & -\frac{\mu'(t)}{\mu(t)^2}\left[ \frac{\a}{\a+1} V +\frac{1}{\a+1}(y \cdot \nabla)V \right]=V_s+(V\cdot \nabla )V +\nabla P,\\
 & \qquad \mathrm{div}\, V=0,\\
 & V(y,0)=V_0(y)=v_0 (y).
  \endaligned \right.
 $$
 We note that the special cases
 $$\mu (t)=\frac{1}{T-t}, \quad \mu (t)=\frac{1}{(T-t)^\g }, \g >1  $$
  are considered in the previous section.
  In this section we choose
$\mu(t)=\exp\left[\pm\g\int_0 ^t \|\nabla v(\tau)\|_{L^\infty}d\tau \right], \g\geq 1.$
Then,
\bq\label{2.4}
 v(x,t)&=&\exp\left[\frac{\pm\g\a}{\a+1}\int_0 ^t\|\nabla v(\tau)\|_{L^\infty}d\tau \right] \ V\left(y, s \right),\\
 \label{2.5}
 p(x,t)&=&\exp\left[\frac{\pm 2\g\a}{\a+1}\int_0 ^t\|\nabla v(\tau)\|_{L^\infty}d\tau \right] \ P\left(y, s\right)
 \eq
 with
\bq \label{2.6}
y&=&\exp\left[\frac{\pm \g}{\a+1}\int_0 ^t\|\nabla v(\tau)\|_{L^\infty}d\tau \right]x, \n \\
 s&=&\int_0 ^t\exp\left[\pm \g\int_0 ^\tau\|\nabla v(\sigma)\|_{L^\infty}d\sigma \right]d\tau
\eq
respectively for the signs $\pm$.
Substituting $(v,p)$ in (\ref{2.4})-(\ref{2.6}) into the $(E)$, we find that $(E_*)$ becomes
 $$
 (E_\pm) \left\{ \aligned
  & \mp \g \|\nabla V(s)\|_{L^\infty}\left[ \frac{\a}{\a+1} V +\frac{1}{\a+1}(y \cdot \nabla)V \right]=V_s+(V\cdot \nabla )V +\nabla P,\\
 & \mathrm{div}\, V=0,\\
 & V(y,0)=V_0(y)=v_0 (y)
  \endaligned \right.
 $$
 respectively for $\pm$.
Similar equations to the system $(E_\pm )$, without the term involving $(y\cdot\nabla)V$ are introduced and studied in \cite{cha4}, where
 similarity type of transform with respect to only time variables was considered.
The argument of the global/local well-posedness of the system $(E_\pm)$ respectively from the local well-posedness result of the Euler equations
is as follows.
  We define
 $$
 S^\pm =\int_0 ^{T}\exp\left[\pm \g\int_0 ^\tau\|\nabla v(\sigma)\|_{L^\infty}d\sigma
 \right]d\tau .
 $$
 Then, $ S^\pm$ is the maximal time of existence of classical solution  for the system $ (E_\pm)$. We also note the following integral invariant of the transform,
 $$
 \int_0 ^{T} \|\nabla v(t)\|_{L^\infty}dt=\int_0 ^{S^\pm} \|\nabla V^\pm (s)\|_{L^\infty}ds.
 $$
  The key advantage of our choice of the function $\mu(t)$ here is that the convection term is dominated
 by $\mp \g \|\nabla V(s)\|_{L^\infty} V$ in the transformed system $(E\pm)$ in the vorticity formulation, which enable us to derive new \textit{a priori} estimates for $\|\o (t)\|_{L^\infty}$ as follows.
\begin{thm}
Given $m>5/2$ and  $v_0 \in H^m (\Bbb R^3)$ with div $v_0=0$, let $\o $ be the vorticity of the solution $v\in C([0, T);H^m (\Bbb R^3 ))$ to the Euler equations (E). Then we have an upper estimate
 \bq\label{2.7}
\|\o (t)\|_{L^\infty}\leq \frac{\|\o_0\|_{L^\infty} \exp\left[ \g \int_0 ^t \|\nabla v (\tau)\|_{L^\infty} d\tau\right]}{1+ (\g-1) \|\o_0 \|_{L^\infty}\int_0 ^t\exp\left[ \g\int_0 ^\tau\|\nabla v(\sigma)\|_{L^\infty}d\sigma \right]d\tau},
 \eq
 and lower one
  \bq\label{2.8}
\|\o (t)\|_{L^\infty}\geq \frac{\|\o_0\|_{L^\infty} \exp\left[ -\g \int_0 ^t \|\nabla v (\tau)\|_{L^\infty} d\tau\right]}{1-(\g-1) \|\o_0 \|_{L^\infty}\int_0 ^t\exp\left[ -\g\int_0 ^\tau\|\nabla v(\sigma)\|_{L^\infty}d\sigma \right]d\tau}
 \eq
  for all $\g \geq 1$ and $t\in [0, T).$  The denominator of the right hand side of (\ref{2.8}) can be  estimated from below as
  \bb\label{2.8a}
  1-(\g-1) \|\o_0 \|_{L^\infty}\int_0 ^t\exp\left[ -\g\int_0 ^\tau\|\nabla v(\sigma)\|_{L^\infty}d\sigma \right]d\tau \geq \frac{1}{(1+\|\o_0 \|_{L^\infty} t)^{\g-1}},
  \ee
  which shows that the finite time blow-up does not follow from (\ref{2.8}).
\end{thm}
{\bf Remark 2.1} We observe that for $\g=1$, the estimates (\ref{2.7})-(\ref{2.8}) reduce to the well-known ones in (\ref{1.11}) with $p=\infty$.
Moreover, combining (\ref{2.7})-(\ref{2.8})  together, we  easily derive another new estimate,
  \bb\label{2.9}
 \frac{\sinh \left[ \g \int_0 ^t \|\nabla v(\tau )\|_{L^\infty} d\tau \right]}{\int_0 ^t \cosh \left[\g \int_\tau ^t \|\nabla v (\sigma )\|_{L^\infty}d\sigma\right] d\tau}\geq (\g-1)\|\o_0 \|_{L^\infty}.
\ee
\ \\
{\bf Proof of Theorem 2.1}
Below we denote $V^\pm$ for the solutions of $(E_\pm)$ respectively, and $\O^\pm = \mathrm{curl}\, V^\pm . $ Note that $V^\pm _0=v_0:=V_0$ and $\O^\pm _0 =\o_0:=\O_0$.
We will first derive the following estimates for the system $(E_\pm)$.
 \bq\label{2.10}
\| \O ^+(s)\|_{L^\infty}&\leq& \frac{\|\O_0\|_{L^\infty}}{1+(\g-1 )s\|\O_0 \|_{L^\infty}},\\
\label{2.11}
\| \O ^-(s)\|_{L^\infty}&\geq &\frac{\|\O_0\|_{L^\infty}}{1-(\g-1 )s\|\O_0 \|_{L^\infty}},
\eq
as long as $V^\pm(s)\in H^m (\Bbb R^3 )$.
Taking curl of the first equation of $(E_\pm)$, we have
\bb\label{2.12}
\mp\g \|\nabla V\|_{L^\infty} \left[\O -\frac{1}{\a +1}(y\cdot \nabla )\O\right]=
\O_s +(V\cdot \nabla )\O -(\O\cdot \nabla )V.
\ee
Multiplying $\Xi=\O/|\O|$ on the both sides of (\ref{2.12}), we deduce
\bq\label{2.13}
\lefteqn{|\O|_s +(V\cdot \nabla )|\O|\mp  \frac{\|\nabla V(s)\|_{L^\infty} }{\a +1}(y\cdot \nabla )|\O|=(\Xi\cdot \nabla V\cdot \Xi \mp\|\nabla V\|_{L^\infty}) |\O|}\hspace{1.5in}\n \\
&&\qquad\mp(\g -1)\|\nabla V\|_{L^\infty} |\O|\n \\
&&\left\{ \aligned &
 \leq -(\g -1)\|\nabla V\|_{L^\infty} |\O|\quad \mbox{for}\,(E_+)\\
 &\geq (\g -1)\|\nabla V\|_{L^\infty} |\O|\quad \mbox{for}\, (E_-),
 \endaligned \right.
\eq
since $|\Xi\cdot \nabla V\cdot \Xi |\leq |\nabla V|\leq \|\nabla V\|_{L^\infty}. $
Given smooth solution $V(y,s)$ of $(E_\pm)$, we  introduce the particle trajectories $\{ Y_\pm (a,s)\}$ defined by
$$
\frac{\partial Y(a,s)}{\partial s}=V_\pm(Y(a,s),s)\mp\frac{\|\nabla V(s)\|_{L^\infty}}{\a+1} Y(a,s) \quad ;\quad
Y(a,0)=a .
$$
Recalling the estimate
$$\|\nabla V(s)\|_{L^\infty}\geq \|\O(s)\|_{L^\infty} \geq |\O (y,s)|\qquad \forall y\in \Bbb R^3,
$$
we can further estimate from (\ref{2.13})
\bb\label{2.14}
\frac{\partial}{\partial s} |\O(Y(a,s),s)| \left\{ \aligned &
 \leq -(\g -1)|\O(Y(a,s),s)|^2 \quad \mbox{for}\, (E_+)\\
 &\geq (\g -1) |\O(Y(a,s),s)|^2\quad \mbox{for}\, (E_-) .
 \endaligned \right.
\ee
Solving these differential inequalities (\ref{2.14}) along the particle trajectories, we obtain that
\bb\label{2.15}
 |\O (Y(a,s) ,s) | \left\{ \aligned &
 \leq \frac{|\O_0 (a)|}{1 +(\g -1)s|\O_0(a )|}\quad \mbox{for}\,(E_+)\\
 &\geq \frac{|\O_0 (a)|}{1 -(\g -1)s|\O_0(a )|}\quad \mbox{for}\, (E_-).
 \endaligned \right.
\ee
Writing the first inequality of (\ref{2.15}) as
$$
 |\O^+ (Y(a,s) ,s) | \leq \frac{1}{\frac{1}{|\O_0 (a)|} +(\g -1)s}
 \leq \frac{1}{\frac{1}{\|\O_0 \|_{L^\infty}} +(\g -1)s},
 $$
 and then taking supremum over $a\in \Bbb R^3$, which is equivalent to taking supremum
 over $Y(a,s)\in \Bbb R^3$ due to the fact that the mapping $a\mapsto Y(a,s)$ is a deffeomorphism(although not volume preserving) on $\Bbb R^3$ as long as $V\in C([0, S); H^m (\Bbb R^3))$, we obtain (\ref{2.10}).
 In order to derive (\ref{2.11})  from the second inequality of (\ref{2.15}), we first write
 $$
 \|\O^-( s)\|_{L^\infty}\geq |\O (Y(a,s),s)|\geq  \frac{1}{\frac{1}{|\O_0 (a) |} -(\g -1)s},
 $$
 and than take supremum over $a\in \Bbb R^3$.
Finally, in order  to obtain  (\ref{2.7})-(\ref{2.8}), we just change variables from (\ref{2.10})-(\ref{2.11})  back to the original physical ones, using the fact
\bqn
\O^+ (y,s)&=&\exp\left[ -\g \int_0 ^t \|\nabla v(\tau )\|_{L^\infty} d\tau\right]\o (x,t),\\
s&=&\int_0 ^t \exp\left[ \g \int_0 ^\tau \|\nabla v(\sigma )\|_{L^\infty} d\sigma\right]d\tau
\eqn
for (\ref{2.7}), while  in order to deduce (\ref{2.8}) from (\ref{2.11}) we substitute
\bqn
\O^- (y,s)&=&\exp\left[ \g \int_0 ^t \|\nabla v(\tau)\|_{L^\infty}d\tau \right]\o (x,t),\\
s&=&\int_0 ^t \exp\left[ -\g \int_0 ^\tau\|\nabla v(\sigma)\|_{L^\infty} d\sigma\right]d\tau .
\eqn
Now we can rewrite (\ref{2.8}) as
$$
\|\o (t)\|_{L^\infty} \geq -\frac{1}{\g-1} \frac{d}{dt} \log\left\{1-(\g-1) \|\o_0 \|_{L^\infty}\int_0 ^t\exp\left[ -\g\int_0 ^\tau\|\nabla v(\sigma)\|_{L^\infty}d\sigma \right]d\tau \right\}.
$$
Thus,
\bq\label{integ}
&&\int_0 ^t \|\nabla v(\tau)\|_{L^\infty} d\tau\geq\int_0 ^t \|\o (\tau)\|_{L^\infty} d\tau\geq \n \\
&&\geq -\frac{1}{\g-1}\log\left\{1-(\g-1) \|\o_0 \|_{L^\infty}\int_0 ^t\exp\left[ -\g\int_0 ^\tau\|\nabla v(\sigma)\|_{L^\infty}d\sigma \right]d\tau \right\}.\n \\
\eq
Setting
$$ y(t):=1-(\g-1) \|\o_0 \|_{L^\infty}\int_0 ^t\exp\left[ -\g\int_0 ^\tau\|\nabla v(\sigma)\|_{L^\infty}d\sigma \right]d\tau,
$$
We find further integrable structure in (\ref{integ}), which is
$$
y'(t)\geq -(\g-1)\|\o_0\|_{L^\infty} y(t)^{\frac{\g}{\g-1}}.
$$
Solving this differential inequality, we obtain (\ref{2.8a}).
$\square$\\
\ \\
{\em In the last part of this section we  fix $\mu(t):=\exp\left[\int_0 ^t \|\nabla v(\tau)\|_{L^\infty}d\tau \right]$.}\\
We assume our local classical solution in $H^m (\Bbb R^3)$ blows up at $T$, and hence
$\mu(T-0)=\exp\left[\int_0 ^T \|\nabla v(\tau)\|_{L^\infty}d\tau \right]=\infty$.
Given $(\a,p)\in (-1, \infty)\times (0, \infty)$,  as previously, we say the blow-up is $\a-$asymptotically self-similar in the sense of $L^p$ if there exists
$\bar{V}=\bar{V}_\a\in \dot{W}^{1,p}(\Bbb R^3 )$ such that the following convergence holds true.
\bb\label{2.16}
\lim_{t\to T}\mu(t)^{-1}\left\|\nabla v(\cdot, t) - \mu(t)\nabla \bar{V}\left(\mu(t)^{\frac{1}{\a+1} }(\cdot)\right)\right\|_{L^\infty}=0
\ee
for $p=\infty$, and
\bb\label{2.17}
\lim_{t\to T}\mu(t)^{-1+\frac{3}{(\a+1)p}}\left\|\o(\cdot ,t)-\mu(t)^{1-\frac{3}{(\a+1)p}} \bar{\O }\left(\mu(t)^{\frac{1}{\a+1}} (\cdot)\right)\right\|_{L^p}=0
\ee
for $p\in (0, \infty)$.
 The above limiting function $\bar{V}$ with $\bar{\O }\neq 0$ is called the blow-up profile as previously.
\begin{pro} Let $\a \neq 3/2$. Then there exists no $\a-$ asymptotically self-similar blow-up in the sense of $L^\infty$ with the blow-up profile belongs to $L^2(\Bbb R^3)$.
\end{pro}
{\bf Proof } Let us suppose that there exists $\bar{V}\in \dot{W}^{1, \infty}(\Bbb R^3)\cap L^2 (\Bbb R^3)$ such that (\ref{2.16}) holds, then we will show that $\bar{V}=0$.
In terms of the self-similar variables (\ref{2.16}) is translated into
$$\lim_{s\to \infty} \|\nabla V(\cdot, s)-\nabla \bar{V}\|_{L^\infty}=0,$$
where $V$ is defined in (\ref{2.1}).
If $\|\nabla \bar{V}\|_{L^\infty}=0$, then, the condition $\bar{V}\in L^2(\Bbb R^3)$ implies that $\bar{V}=0$, and there is noting to prove. Let us suppose $\|\nabla \bar{V}\|_{L^\infty}>0$.
The equations satisfied $\bar{V}$ are
\bb\label{2.18}
\left\{ \aligned
&-\|\nabla \bar{V}\|_{L^\infty}\left[ \frac{\a}{\a+1} \bar{V }+\frac{1}{\a+1}(y \cdot \nabla)\bar{V} \right]=(\bar{V}\cdot \nabla )\bar{V }+\nabla \bar{P},\\
&\qquad\mathrm{div} \,\bar{ V}=0
\endaligned
\right.
\ee
for a scalar function $\bar{P}$. Taking $L^2 (\Bbb R^3)$ inner product of the first equation of (\ref{2.18}) by $\bar{V}$ we obtain
$$
\frac{\|\nabla \bar{V}\|_{L^\infty}}{\a+1} \left(\a-\frac{3}{2}\right) \|\bar{V}\|_{L^2}=0.
$$
Since $\|\nabla \bar{V}\|_{L^\infty}\neq 0$ and $\a \neq \frac 32$, we have $\|\bar{V}\|_{L^2}=0$, and $\bar{V}=0$. $\square$\\
\ \\
 \begin{pro}
There exists no $\a-$asymptotically self-similar blowing up solution to (E) in the sense of $L^p$ if
$0<p < \frac{3}{2(\a+1)}.$
 \end{pro}
 {\bf Proof }
 Suppose there exists $\a-$asymptotically self-similar blow-up at $T$ in the sense of $L^p$. Then, there exists $\bar{\O} \in L^p (\Bbb R^3)$ such that, in terms of the self-similar variables introduced in (\ref{2.1})-(\ref{2.2}), we have
 \bb\label{2.19}
 \lim_{s\to \infty} \|\O (s)\|_{L^p}=\|\bar{\O}\|_{L^p} <\infty.
 \ee
   We represent the $L^p$ norm of $\|\o (t)\|_{L^p}$ in terms of similarity variables to obtain
\bb\label{2.20}
\|\o (t)\|_{L^p}
= \mu(t)^{1-\frac{3}{(\a+1)p}}\|\O(s)\|_{L^p}, \quad \mu (t)=\exp\left[ \int_0 ^t \|\nabla v (\tau)\|_{L^\infty} d\tau\right].
\ee
Substituting this into the lower estimate part of (\ref{1.11}), we have
 \bb\label{2.21}
\mu(t)^{-2+ \frac{3}{(\a+1)p}}\leq \frac{\|\O (s)\|_{L^p}}{\|\O_0 \|_{L^p}}.
 \ee
 If $-2+\frac{3}{(\a+1)p}>0$, then taking $t\to T$ the above inequality we obtain,
\bqn
 \lefteqn{\infty=\lim\sup_{t\to T}\mu(t)^{-2+\frac{3}{(\a+1)p}}\|\O_0 \|_{L^p}}\hspace{.0in}\n \\
&&\leq\lim\sup_{s\to \infty}\|\O (s)\|_{L^p}=\|\bar{\O}\|_{L^p},
\eqn
which is a contradiction to (\ref{2.19}). $\square$\\
\ \\
\section{\bf The case of the 3D Navier-Stokes equations}
 \setcounter{equation}{0}

In this section we concentrate on the following 3D Navier-Stokes equations in $\Bbb R^3$ without forcing term.
\[
\mathrm{ (NS)}
 \left\{ \aligned
 &\frac{\partial v}{\partial t} +(v\cdot \nabla )v =\Delta  v-\nabla p ,
 \quad (x,t)\in {\Bbb R^3}\times (0, \infty) \\
 &\quad \textrm{div }\, v =0 , \quad (x,t)\in {\Bbb R^3}\times (0,
 \infty)\\
  &v(x,0)=v_0 (x) \quad x\in \Bbb R^3.
  \endaligned
  \right.
  \]
First, we exclude asymptotically self-similar
   singularity  of type II of (NS), for which  the blow-up rate is given by (\ref{type2}).
  We have the following theorem.
  \begin{thm}
Let $p\in [3, \infty)$ and $v\in C([0, T); L^p (\Bbb R^3))$ be a local classical solution of the Navier-Stokes equations constructed by Kato(\cite{kat1}). Suppose  there exists  $\g >1$ and $\bar{V}\in L^p (\Bbb R^3)$ such that the following convergence holds true.
\bb\label{3.1}
\lim_{t\to T}(T-t)^{\frac{(p-3)\g }{2p}}\left\|v(\cdot, t)-(T-t)^{-\frac{(p-3)\g }{2p}}\bar{V}\left(\frac{\cdot}{(T-t)^{\frac{\g}{2}}} \right) \right\|_{L^p}=0,
\ee
If the blow-up profile $\bar{V}$ belongs to $\dot{H}^1 (\Bbb R^3)$, then  $\bar{V}=0$.
\end{thm}
{\bf Proof } Since the main part of the proof is essentially identical to that of Theorem 1.3, we will be brief. Introducing the self-similar variables of the form (\ref{1.17})-(\ref{1.19}) with $\a=\frac12$, and substituting $(v,p)$ into the Navier-Stokes equations, we find that $(V,P)$ satisfies
$$
  \left\{ \aligned
  & -\frac{\g}{2s(\g-1) +2T^{1-\g}}\left[  V +(y \cdot \nabla)V \right]=V_s+(V\cdot \nabla )V -\Delta V+\nabla P,\\
 & \quad \mathrm{div}\, V=0,\\
 & V(y,0)=V_0(y)=v_0 (y).
  \endaligned \right.
$$
The hypothesis (\ref{3.1}) is now translated  as
$$
\lim_{s\to \infty}\|V(\cdot, s)-\bar{V}(\cdot )\|_{L^p}=0
$$
 Following exactly same argument as in the proof of Theorem 1.3, we can deduce that $\bar{V}$ is a stationary solution of the Navier-Stokes equations, namely there exists $\bar{P}$ such that
\bb\label{3.1a}
(\bar{V}\cdot \nabla )\bar{V}=\Delta \bar{V }-\nabla \bar{P}, \qquad \mathrm{div }\, \bar{V}=0.
\ee
In the case $\bar{V}\in \dot{H}^1\cap L^p (\Bbb R^3)$, we easily from (\ref{3.1a}) that
$\int_{\Bbb R^3} |\nabla \bar{V}|^2 dy=0$, which implies  $\bar{V}=0$. $\square$\\
 \ \\
 Next, we derive a new a priori estimates for classical solutions of the 3D Navier-stokes equations.
 \begin{thm}
Given $v_0 \in H^1(\Bbb R^3)$ with div $v_0=0$, let $\o $ be the vorticity of the classical solution $v\in C([0, T);H^1 (\Bbb R^3 ))\cap C((0, T); C^\infty (\Bbb R^3))$ to the  Navier-Stokes  equations (NS). Then,  there exists an
absolute constant $C_0>1$ such that for all $\g \geq C_0$ the following enstrophy estimate holds true.
 \bq\label{3.2}
\|\o (t)\|_{L^2}\leq \frac{\|\o_0\|_{L^2} \exp\left[ \frac{\g}{4} \int_0 ^t \|\o (\tau)\|_{L^2}^4 d\tau\right]}{\left\{ 1+ (\g-C_0) \|\o_0 \|_{L^2}^4\int_0 ^t\exp\left[ \g\int_0 ^\tau\|\o (\sigma)\|_{L^2}^4d\sigma \right]d\tau\right\}^{\frac14}}.
 \eq
The denominator of (\ref{3.2}) is estimated from below by
 \bb\label{3.2a}
 1+ (\g-C_0) \|\o_0 \|_{L^2}^4\int_0 ^t\exp\left[ \g\int_0 ^\tau\|\o (\sigma)\|_{L^2}^4d\sigma \right]d\tau
 \leq \frac{1}{(1-C_0 \|\o_0 \|_{L^2}^4 t )^{\frac{\g-C_0}{C_0}}}
 \ee
 for all $\g \geq C_0$.
 \end{thm}
 \noindent{\bf Proof } Let $(v,p)$ be a classical solution of the Navier-Stokes equations, and $\o$ be its vorticity.
 We transform from $(v,p)$ to $(V, P)$ according to the formula, given by (\ref{2.1})-(\ref{2.3}), where
 $$
 \mu (t)=\exp\left[ \g\int_0 ^t \|\o (\tau )\|_{L^2} ^4 d\tau\right].
 $$
Substituting (\ref{2.1})-(\ref{2.3}) with such $\mu(t)$ into (NS), we obtain the equivalent equations satisfied by $(V,P)$
  $$
 (NS_*) \left\{ \aligned
  & \frac{-\g\|\O (s)\|_{L^2}^4}{2}\left[  V +(y \cdot \nabla)V \right]=V_s+(V\cdot \nabla )V -\Delta V-\nabla P,\\
 & \qquad \mathrm{div}\, V=0,\\
 & V(y,0)=V_0(y)=v_0 (y).
  \endaligned \right.
 $$
Operating curl on the evolution equations of $(NS_*)$, we obtain
 \bb\label{3.5}
 \frac{-\g \|\O (s)\|_{L^2}^4}{2}\left[  2\O +(y \cdot \nabla)\O\right]=\O_s+(V\cdot \nabla )\O-(\O \cdot \nabla )V-\Delta \O.
 \ee
 Taking $L^2(\Bbb R^3)$ inner product of (\ref{3.5}) by $\O$, and integrating by part, we
 estimate
 \bq\label{3.6}
 \lefteqn{\frac12 \frac{d}{ds} \|\O\|_{L^2}^2 +\|\nabla \O \|_{L^2}^2 + \frac{\g}{4} \|\O\|_{L^2}^6
 =\int_{\Bbb R^3} (\O \cdot \nabla )V \cdot \O dy}\hspace{1.in}\n \\
 &&\leq\|\O\|_{L^3} \|\nabla V\|_{L^2} \|\O\|_{L^6}\leq C \|\O \|_{L^2}^{\frac32} \|\nabla \O\|_{L^2}^{\frac32}\n \\
 &&\leq  \|\nabla \O \|_{L^2}^2 +\frac{C_0}{4} \|\O\|_{L^2}^6\n \\
 \eq
 for an absolute constant $C_0>1$, where we used the fact $\|\O\|_{L^2}=\|\nabla V\|_{L^2}$, the Sobolev imbedding, $\dot{H}^1(\Bbb R^3)\hookrightarrow L^6(\Bbb R^3)$, the Gagliardo-Nirenberg inequality in $\Bbb R^3$,
 $$ \|f\|_{L^3}\leq C \|f\|_{L^2} ^{\frac12}\|\nabla f\|_{L^2}^{\frac12}.$$
 and Young's inequality of the form $ab\leq a^p/p+b^q/q,$ $1/p+1/q=1$.
 Absorbing  the term $\|\nabla \O \|_{L^2}^2$ to the left hand side,  we have from (\ref{3.6})
 \bb\label{3.7}
\frac{d}{ds} \|\O\|_{L^2}^2\leq -\frac{\g -C_0}{2} \|\O\|_{L^2}^6.
\ee
Solving the differential inequality (\ref{3.7}), we have
\bb\label{3.8}
\|\O (s)\|_{L^2}\leq \frac{\|\O_0 \|_{L^2} }{\left[ 1+(\g-C_0) s \|\O_0\|_{L^2}^4\right]^{\frac14}}.
\ee
Transforming back to the original variables and functions, using the relations
\bqn
s&=&\int_0 ^t \exp\left[ \g \int_0 ^\tau \|\o (\sigma )\|_{L^2} ^4 d\sigma\right] d\tau,\\
\|\o (t)\|_{L^2}&=& \|\O (s)\|_{L^2} \exp\left[ \frac{\g}{4} \int_0 ^t \|\o (\tau)\|_{L^2}^4 d\tau\right],
\eqn
we obtain (\ref{3.2}).
Next, we observe (\ref{3.2}) can be written as
 $$\|\o (t)\|_{L^2}^4\leq \frac{1}{(\g-C_0 )} \frac{d}{dt} \log\left\{ 1+ (\g-C_0) \|\o_0 \|_{L^2}^4\int_0 ^t\exp\left[ \g\int_0 ^\tau\|\o (\sigma)\|_{L^2}^4d\sigma \right]d\tau\right\},
 $$
 which, after integration over $[0,t]$, leads to
 \bq\label{cla1}
\int_0 ^t \|\o (\tau )\|_{L^2}^4 d\tau
 \leq \frac{1}{(\g-C_0 )}\log \left\{ 1+ (\g-C_0) \|\o_0 \|_{L^2}^4\int_0 ^t\exp\left[ \g\int_0 ^\tau\|\o (\sigma)\|_{L^2}^4d\sigma \right]d\tau\right\}\n \\
 \eq
 for all $\g >C_0$. Setting
 $$y(t):=1+ (\g-C_0) \|\o_0 \|_{L^2}^4\int_0 ^t\exp\left[ \g\int_0 ^\tau\|\o (\sigma)\|_{L^2}^4d\sigma \right]d\tau,
 $$
 we find that (\ref{cla1}) can be written in the form of a  differential inequality,
 $$ y'(t)\leq (\g-C_0 )\|\o_0\|_{L^2}^4 \,y(t)^{\frac{\g}{\g-C_0}},
 $$
 which can be solved to provide us with (\ref{3.2a}).
 $\square$

\end{document}